\documentclass[11pt]{article}
\newtheorem{df}{Definition}[section]
\newtheorem{theo}[df]{Theorem}
\newtheorem{prop}[df]{Proposition}
\newtheorem{cor}[df]{Corollary}
\newtheorem{lem}[df]{Lemma}
\newtheorem{rem}[df]{Remark}
\newenvironment{proof}{{\sc Proof\ }}{\hfill $\Box$}
\usepackage{latexsym,amssymb}
\usepackage{amsmath}

\begin{document}

\title{A classification of\\ spherical symmetric
$CR$ manifolds}
\author{G. Dileo and A. Lotta}

\date{}

\maketitle
\begin{abstract}
In this paper we classify the simply connected, spherical
pseudohermitian manifolds whose Webster metric is $CR$-symmetric.
\end{abstract}

\medskip
\noindent
{\small

}

\section{Introduction}

A spherical $CR$ manifold is a strongly pseudoconvex $CR$ manifold
$(M,HM,J)$ of hypersurface type which is locally $CR$-equivalent
to the sphere $\mathbb S^{2n+1}$, $n=\dim_{CR}M$, endowed with the
standard $CR$ structure as a real hypersurface of $\mathbb
C^{n+1}$. Recall that strong pseudoconvexity means positive
definiteness of the Levi form $L_\eta$ associated to a suitable
global section $\eta$ of the annihilator $H^oM$ of the holomorphic
tangent bundle of $M$. The 1-form $\eta$ is usually called a {\em
pseudohermitian structure} on $M$, and it canonically determines a
Riemannian metric $g_\eta$ which is compatible with the partial
complex structure $J:HM\to HM$ (cf. e.g. \cite{Web0},\cite{TAN}).
We shall call $g_\eta$ the {\em Webster metric} associated to
$\eta$. Denoting by $\xi$ the Reeb vector field of the contact
form $\eta$, at each point $x\in M$ we have an orthogonal
decomposition
$$T_xM=H_xM\oplus\mathbb R\xi_x$$
with respect to $g_\eta$, and moreover
${g_\eta}_{|H_xM}=(L_\eta)_x$, $\,\,g_\eta(\xi,\xi)=1$.

Spherical $CR$ manifolds are characterized by $B=0$,
where $B$ is the Chern-Moser-Tanaka pseudoconformal
invariant tensor field of type $(1,3)$, and they
represent flat spaces among strongly pseudoconvex
$CR$ manifolds from the point of view
of Cartan geometry (see e.g. \cite{tanhyp}).
The simply connected, homogeneous spherical hypersurfaces
of the Euclidean space $\mathbb C^{n+1}$ were fully
classified by D. Burns and S. Shnider in \cite{BS}.
In particular, it is known that the unique
compact simply connected homogeneous spherical
hypersurface of $\mathbb C^{n+1}$ is $\mathbb S^{2n+1}$,
up to $CR$-equivalence.
In this paper we adopt a geometric point of view
in studying spherical $CR$ manifolds, concentrating
our attention to $CR$-{\em symmetric} Webster metrics
$g_\eta$. For the general notion of a symmetric Hermitian
metric on a $CR$ manifold we refer to \cite{KZ}.
Here we recall that a Webster metric $g_\eta$
is $CR$-symmetric if  for each point $x\in M$
there exists a  $CR$-isometry
$\sigma:M\to M $
with $\sigma(x)=x$ and
$${(d\sigma)_x}_{|H_xM}=-Id.$$
See \S\ref{CR-sym} for more details.

Actually the standard metric $g_o$ of curvature
1 on the sphere $\mathbb S^{2n+1}$ is a $CR$-symmetric
Webster metric $g_{\eta_o}$ for the choice
of a canonical contact form $\eta_o$.  The
symmetry at a point $x\in\mathbb S^{2n+1}$ is the restriction
of the unitary reflection
$\sigma_x(z)=2<x,z>x-z$ with respect to the
standard Hermitian scalar product of
$\mathbb C^{n+1}$ (cf. \cite{KZ}).

More generally, any {\em Sasakian space form}
(see e.g. \cite{BLAIR+})
is a spherical $CR$-symmetric pseudohermitian
manifold.
Indeed, in the Sasakian case, the Webster metric $g$
is $CR$-symmetric if and only if $M$ is
a $\varphi$-symmetric space (for this notion see
e.g. \cite{TA}).
From the classification of Sasakian
$\varphi$-symmetric spaces
carried out by J.A. Jim\'enez and O. Kowalski in \cite{JK}
we also see that for $n\ge 2$ and $0<k<n$
there exists a principal fiber
bundle $P^n_k\to{\mathbb C\mathbb P}^k\times{\mathbb C\mathbb H}^{n-k}$
with $\varphi$-symmetric Sasakian total space,
the base space $N={\mathbb C\mathbb P}^k\times{\mathbb C\mathbb H}^{n-k}$
being the product of two K\"ahler space forms
with holomorphic curvatures $1$ and $-1$ respectively.
This $CR$-symmetric space $P^n_k$
is spherical since the base manifold
is Bochner-flat (cf. \cite{Bryant}), according to a result of
S. Webster which identifies the Chern-Moser tensor
of $P^n_k$ with the Bochner tensor of $N$
(see \cite{Web} or \cite{David}).

\smallskip
In this paper we get a complete classification
of the simply connected spherical, $CR$-symmetric
pseudohermitian manifolds.
We say that two pseudohermitian manifolds
$(M_1,HM_1,J_1,\eta_1)$ and
$(M_2,HM_2,J_1,\eta_2)$ are {\em  homothetic}
if there exists a $CR$-diffeomorphism $f:M_1\to M_2$
such that $f^*\eta_2=\alpha\eta_1$ with $\alpha$
a positive constant.

\medskip
Our main result is the following:

\bigskip
\begin{theo}
\label{main} Every
simply connected, spherical
$CR$-symmetric
pseudohermitian manifold  of $CR$-dimension $n\ge 2$
is homothetic to one of the following spaces:
$$\Bbb S^{2n+1},\,
H^{2n+1},
B^n\times\Bbb R,
T_1\Bbb H^{n+1},\,
P^n_1,\dots,P^n_{n-1}.
$$
\end{theo}

The first three spaces are the
simply connected Sasakian space forms
as described in \cite{BLAIR+}, p. 114.
$H^{2n+1}$ denotes the Heisenberg group endowed
with its standard Webster flat Sasakian structure, while
$B^n\times\Bbb R$ is the product of a
K\"ahler bounded domain in $\mathbb C^n$ having
constant holomorphic negative curvature with the real line,
which carries a Sasakian structure
with constant $\varphi$-sectional curvature $<-3$.
The fourth space is the tangent sphere bundle $T_1{\mathbb H}^{n+1}$ of the
Riemannian space form of curvature $-1$, with
its standard $CR$ structure and pseudohermitian
structure studied for example in \cite{TanBundle}.
The remaining $n-1$ spaces
are the Sasakian $\varphi$-symmetric spaces described above.

\smallskip
As a consequence we obtain

\begin{cor}
\label{cor1}
Up to homothety, the sphere $\Bbb S^{2n+1}$
is the unique simply connected, compact, spherical
$CR$-symmetric pseudohermitian manifold having $CR$-dimension
$n\ge 2$.
\end{cor}

\begin{cor}
\label{cor2}
A spherical, $CR$-symmetric pseudohermitian
manifold of $CR$-dimen\-sion
$n\ge 2$, having positive pseudoholomorphic
curvature at some point, is compact and is actually a Sasakian pseudohermitian
space form.
\end{cor}

The notion of pseudoholomorphic sectional curvature
is described in detail in \S\ref{prelim}.

\smallskip
The key facts for establishing Theorem \ref{main} are provided
by the following result,
where basic features
of $CR$-symmetric Webster metrics are analyzed
(see Theorems \ref{equivalent} and \ref{B=0}):

\begin{theo}
\label{kmu}
Let $(M,HM,J,\eta)$ be a pseudohermitian
manifold of $CR$ dimension $n\ge 2$.
Assume that the Webster metric $g:=g_\eta$ is
not Sasakian. Then

\medskip\noindent
a) $(M,HM,J,g)$ is locally $CR$-symmetric 
if and only if the underlying contact metric structure
$(\varphi,\xi,\eta,g)$ satisfies
the $(k,\mu)$-nullity condition, that is
$$R(X,Y)\xi=k(\eta(Y)X-\eta(X)Y)+\mu(\eta(Y)hX-\eta(X)hY)\quad\,\,
k,\mu\in\mathbb R$$
where $R$ is the curvature tensor of $g$
and $h=\frac{1}{2}{\cal L}_\xi\varphi$.

\medskip\noindent
b) Assume that $g$ is locally $CR$-symmetric.
Then the following are equivalent:

\smallskip
i)  $M$ is spherical.

ii) The Webster scalar curvature vanishes.

iii) The pseudoholomorphic sectional curvature $\tilde K$
is constant.

iv) $\mu=2$.

\noindent
When one of these equivalent conditions holds,
then $\tilde K=0$ but $\tilde R\not=0$.
\end{theo}

This result provides an interpretation in
$CR$ geometry of the theory of contact Riemannian $(k,\mu)$-spaces
started in \cite{BKP} and fully developed in the
last decade by several authors, especially
by E. Boeckx in \cite{BO1}, \cite{BO2}, \cite{BO3}.
We also remark that Theorem \ref{main} includes the
classification of pseudo-parallel strongly pseudoconvex
$CR$ manifolds with constant pseudoholomorphic
sectional curvature obtained by J.T. Cho in \cite{Cho}.
These manifolds are exactly
$\Bbb S^{2n+1},\,
H^{2n+1},
B^n\times\Bbb R$ and $T_1{\mathbb H}^{n+1}$.
The reason stands in the fact that {\em any} pseudohermitian manifold
with constant pseudoholomorphic
sectional curvature {\em must}  be spherical.
We prove this in \S \ref{sectboch} (Theorem \ref{space form}).

As an application of Theorem \ref{kmu}, in the
last section we study the $CR$ geometry of tangent sphere bundles
of arbitrary constant radius over Riemannian
manifolds with constant curvature. Our approach is
slightly more general than the one appearing in \cite{Tash}
and \cite{TanBundle}.
We show that if $M$ is a hyperbolic Riemannian space form,
each $T_rM$ carries a one parameter family of $CR$-symmetric
non homothetic pseudohermitian structures, exactly
one of which is spherical (Theorem \ref{lambdar}).
This should compared with
the relevant fact that a homogeneous $CR$ manifold which is
{\em homeomorphic} to a sphere, admits a {\em unique} homogeneous
$CR$ structure
\cite{FL}.
Our examples emphasize
that ``homeomorphic'' cannot be replaced by ``homotopically equivalent'',
even if the homogeneous $CR$ structure is spherical.


\smallskip
We also obtain the following:

\begin{cor}
\label{K><0}
Let $(M,g)$ be a Riemannian manifold of constant curvature $K$
and dimension $n\ge 3$.
Consider the standard almost complex structure $J$ on $TM$
defined by
$$JX^H=X^V,\,JX^V=-X^H\quad X\in\frak X(M)$$
where $X^H$ and $X^V$ denote horizontal and vertical lifts.
Endow each $T_rM$, $r>0$, with the induced $CR$ structure
and the standard pseudohermitian structure $\eta_r$.
Then

\smallskip
a) $K<0$ if and only if there exists $r>0$ such that $T_rM$ is
a spherical $CR$ manifold.

\smallskip
b) $K>0$ if and only if there exists $r>0$ such that the
Webster metric $g_{\eta_r}$ is
Sasakian.

\smallskip
c) $K=0$ if and only if  $(T_rM,H(T_rM),J,\eta_r)$ and
$(T_{r'}M,H(T_{r'}M),J,{\eta_{r'}})$ are locally homothetic
pseudohermitian manifolds for each $r,r'>0$.

\medskip
Moreover, if one of the equivalent conditions in a) holds, there exists
a unique $r_o$ such that
$T_{r_o}M$ is spherical, which is related to $K$ by
$$K=-\frac{1}{r_o^2}.$$
When one of the equivalent conditions in b) holds there exists
a unique $r_o$ such that
$g_{\eta_{r_o}}$ is Sasakian, which is related to $K$ by
$$K=\frac{1}{r_o^2}.$$
\end{cor}

\section{Preliminaries}
\label{prelim}
Let $M$ be a connected $\cal{C}^\infty$ manifold of
dimension $2n+k$, $n,k\geq 1$.
A {\em partial complex structure}
of $CR$-dimension $n$
and $CR$-codimension $k$
is a pair $(HM, J)$ where $HM$ is a
smooth real subbundle of the tangent bundle $TM$ having
rank $2n$, and $J$ is a
smooth bundle isomorphism $J:HM\rightarrow HM$, such that $J^2=-I$.
An \emph{almost $CR$ structure} is a partial complex structure
such that
\begin{equation}\label{CR1}[X,Y]-[JX,JY]\in\mathcal{D}\end{equation}
for any $X,Y\in\mathcal{D}$,
where $\mathcal{D}$ denotes the module of all
smooth sections of $HM$.
If, in addition, the formal integrability condition
\begin{equation}\label{CR2}[JX,JY]-[X,Y]-J([JX,Y]+[X,JY])=0\end{equation}
is satisfied, $(M,HM,J)$ is termed
a \emph{$CR$ manifold}.
In this paper we shall be concerned only with the case where
$M$ has $CR$-codimension 1 (hypersurface type).
Assuming (\ref{CR1}) and that
$M$ is orientable, it is known that there exist globally defined nowhere zero
$1$-forms $\eta$ such that $\mathrm{Ker}(\eta)=HM$. The
corresponding {\em Levi form}  is defined by
\[L_\eta(X,Y)=-d\eta(X,JY),\quad X,Y\in{\cal D}.\] The
almost $CR$ structure is said to be \emph{nondegenerate} if $L_\eta$ is
nondegenerate for some $\eta$. In this case, the $1$-form $\eta$ is
a contact form, in the sense that $\eta\wedge(d\eta)^n$ is a volume
form on $M$. Moreover, there exists a unique nowhere vanishing globally
defined vector field $\xi$ transverse to $HM$ such that
\[\eta(\xi)=1,\quad d\eta(\xi,X)=0\]
for any $X\in {\frak X}(M)$. The second condition is equivalent to
$[\xi,{\cal D}]\subset{\cal D}$ or ${\cal L}_\xi\eta=0,$ where
${\cal L}_\xi$ denotes the Lie differentiation with respect to
$\xi$.

An almost $CR$ structure is said to be \emph{strongly pseudoconvex} if
$L_{\eta}$ is positive definite for some $\eta$.
In this case the Levi
form can be canonically extended to a Riemannian metric on $M$, called the
\emph{Webster metric}, defined by
\[g_\eta(X,Y)=L_\eta(X,Y),\quad g_\eta(X,\xi)=0,\quad g_\eta(\xi,\xi)=1,\]
for any $X,Y\in{\cal D}$.  Such a 1-form $\eta$ will be called
a {\em pseudohermitian structure}.

By a {\em pseudohermitian manifold} we shall mean a
strongly pseudoconvex $CR$ manifold $(M,HM,J,\eta)$ on which
a pseudohermitian structure has been fixed.
The partial complex structure $J$ of a pseudohermitian
manifold can be canonically extended
to a tensor
field $\varphi$ of type $(1,1)$ on $M$ such that $\varphi(\xi)=0$
and $\varphi X=JX$ for any $X\in{\cal D}$, which is
an $f$-structure with rank $2n$.
The tensors $(\varphi,\xi,\eta,g_\eta)$ make up a
{\em contact metric structure} on
$M$ in the sense of \cite{BLAIR+}.
Conversely, if $M$ is a contact metric manifold with structure
$(\varphi,\xi,\eta,g)$, then $M$ admits a strongly
pseudoconvex almost $CR$ structure
given by $HM=\mathrm{Im}(\varphi)$ and $J=\varphi|_{HM}$.
The Webster metric $g_\eta$ coincides with $g$.
S. Tanno
proved that this almost $CR$ structure is a $CR$ structure
if and only if
\[(\nabla_X\varphi)Y=g(X+hX,Y)\xi-\eta(Y)(X+hX)\] for any
$X,Y\in{\frak X}(M)$ (cf. \cite{TAN}).
Here $\nabla$ is the Levi-Civita connection of the Webster
metric $g$ and $h$ is the symmetric operator
$h:=\frac{1}{2}\mathcal{L}_\xi\varphi$.

We also recall that a pseudohermitian manifold $(M,HM,J,\eta)$ such
that $h=0$ is called a {\em Sasakian manifold} and the metric
$g_\eta$ is called a {\em Sasakian metric}.

\medskip
Next we recall a special class of contact metric manifolds
with which we will be concerned in the following,
the so-called {\em $(k,\mu)$-spaces}.
Such a space is characterized by
the following property of the Riemannian curvature tensor,
known in the literature as the {\em $(k,\mu)$-nullity condition}:
\begin{equation}\label{(k,mu)}
R(X,Y)\xi=k(\eta(Y)X-\eta(X)Y)+\mu(\eta(Y)hX-\eta(X)hY)\quad
k,\,\mu\in\mathbb R
\end{equation}
for any $X,Y\in{\frak X}(M)$.
In \cite{BKP} the authors prove the relation $h^2=(k-1)\varphi^2$
which implies that $k\leq 1$. If
$k=1$, then $h=0$ and $M$ is a Sasakian manifold. If $k<1$, the
contact metric structure is not Sasakian and $M$ admits three
mutually orthogonal integrable distributions $D(0)=\mathbb{R}\xi,$
$D(\lambda)$ and $D(-\lambda)$, determined by the eigenspaces of
$h$, with $\lambda=\sqrt{1-k}$. Moreover, when $k<1$ the curvature
tensor is completely determined by the condition (\ref{(k,mu)}),
and its explicit expression is the following \cite{BO1}:
\begin{eqnarray}
\label{Rkmu}
R(X,Y)Z &=& \left(1-\frac{\mu}{2}\right)(g(Y,Z)X-g(X,Z)Y)\\
&&{}+g(Y,Z)hX-g(X,Z)hY+g(hY,Z)X-g(hX,Z)Y\nonumber\\
&&{}+\frac{1-\frac{\mu}{2}}{1-k}(g(hY,Z)hX-g(hX,Z)hY)\nonumber\\&&{}-
\frac{\mu}{2}(g(\varphi Y,Z)\varphi X-g(\varphi X,Z)\varphi Y)
+\mu g(\varphi X,Y)\varphi Z\nonumber\\
&&{}+\frac{k-\frac{\mu}{2}}{1-k}(g(\varphi hY,Z)\varphi hX-g(\varphi
hX,Z)\varphi hY)\nonumber\\
&&+\eta(X)\left((k-1+\frac{\mu}{2})g(Y,Z)+(\mu-1)g(hY,Z)\right)\xi\nonumber\\
&&-\eta(Y)\left((k-1+\frac{\mu}{2})g(X,Z)+(\mu-1)g(hX,Z)\right)\xi\nonumber\\
&&-\eta(X)\eta(Z)\left((k-1+\frac{\mu}{2})Y+(\mu-1)hY\right)\nonumber\\
&&+\eta(Y)\eta(Z)\left((k-1+\frac{\mu}{2})X+(\mu-1)hX\right).\nonumber
\end{eqnarray}

In \cite{BO2}, E. Boeckx introduced the invariant
$$I=\frac{1-\mu/2}{\sqrt{1-k}}$$
and proved that two non Sasakian $(k,\mu)$-spaces are locally
homothetic pseudohermitian manifolds if and only if their
invariants coincide.

Moreover, the non Sasakian $(k,\mu)$-spaces are also characterized
by the requirement that
\begin{equation}
\label{etap}
g((\nabla_Xh)Y,Z)=0
\end{equation}
for any $X,Y,Z\in\cal D$. This fact has been proved in \cite{BC},
where contact metric structures satisfying (\ref{etap}) are called
{\em $\eta$-parallel}.

We end this section recalling some basic facts
about the Tanaka-Webster connection. We shall refer
to \cite{SakTak1}.

\begin{theo}\label{connection}
Let $(M,HM,J,\eta)$ be a pseudohermitian manifold
with subordinate contact metric structure
$(\varphi,\xi,\eta,g)$. There is a unique linear connection
$\tilde\nabla$ such that
\begin{equation}\label{canonical}
\tilde\nabla\varphi=0,\quad \tilde\nabla\xi=0,\quad
\tilde\nabla\eta=0,\quad \tilde\nabla g=0,\quad \tilde T_{\cal
D}=0,\quad F=-\frac{1}{2}\varphi {\cal L}_\xi\varphi,
\end{equation}
where $\tilde T$ is the torsion tensor field of
$\tilde\nabla$, $\tilde T_{\cal D}(X,Y)$ denotes the $\cal D$-component
of $\tilde T(X,Y)$ for any $X,Y\in{\cal D}$, and $F$ is the tensor
field of type $(1,1)$ defined by $FX=\tilde T(\xi,X)$ for any
$X\in{\frak X}(M)$.
\end{theo}

The linear connection in the above statement is called the
\emph{canonical connection} or the \emph{Tanaka-Webster connection}
of the pseudohermitian manifold $M$.
Denoting by
$\nabla$ the Levi-Civita connection of $g$, then
$\tilde\nabla=\nabla+H$, with
\begin{equation}\label{H}
H(X,Y)=g(X,\varphi Y)\xi+\eta(X)\varphi Y+\eta(Y)\varphi X
+g(FX,Y)\xi-\eta(Y)FX
\end{equation}
for any $X,Y\in{\frak X}(M)$. The curvature tensor $\tilde R$ of
$\tilde\nabla$ satisfies
\[\tilde R(X,Y)\xi=0,\quad \tilde R(X,Y)\varphi=\varphi\tilde R(X,Y),\quad
\tilde R(X,Y){\cal D}\subset{\cal D}\] for any $X,Y\in{\frak X}(M)$.

If $\sigma\subset H_xM$ is a holomorphic 2-plane in $x\in M$, that
is $J\sigma=\sigma$, the quantity
\[\tilde K (\sigma)= \tilde R_x(X,JX,X,JX)\]
where $\{X,JX\}$ is an orthonormal basis of $\sigma$,
depends only on $\sigma$ and
will be called the \emph{pseudoholomorphic sectional
curvature} of $\sigma$. If $\tilde K (\sigma)$ does not depend
on $\sigma$ and on the point $x$, $M$ will be called a
{\em pseudohermitian space form}.

\section{$CR$-symmetric  Webster metrics}
\label{CR-sym}

Let $(M,HM,J,g)$ be a {\em Hermitian} almost $CR$ manifold, i.e. an
almost $CR$ manifold, having $CR$-codimension $k\ge 1$,
on which a Riemannian metric $g$ is fixed,
whose restriction to $HM$ is Hermitian with respect to $J$.
Denote by ${\cal D}_\infty\subset\frak X(M)$ the Lie
algebra generated by ${\cal D}$.
Let $\sigma:M\rightarrow M$ be an isometric $CR$-diffeomorphism. Then
$\sigma$ is called a \textit{symmetry} at the point $x\in M$ if
$x$ is a fixed point of $\sigma$ and the differential of $\sigma$
at $x$ coincides with $-Id$  on the subspace
${\cal D}_\infty(x)^\perp\oplus H_xM$ of $T_xM$. Here
${\cal D}_\infty(x)=\{X_x|X\in{\cal D}_\infty\}$.

A connected Hermitian almost $CR$ manifold
$M$ is called a (globally) {\em $CR$-symmet\-ric space}
if for each point $x\in M$ there exists a symmetry $\sigma_x$ at $x$
(cf. \cite{KZ}). We shall also say that $g$ is a $CR$-symmetric
Hermitian metric on $(M,HM,J)$.
Since the symmetry at $x$ in uniquely determined (cf. Theorem 3.3
in \cite{KZ}) it makes
sense also to define {\em locally} $CR$-symmetric spaces in
a natural manner.
Observe that, since the symmetries are $CR$ maps,
the integrability condition (\ref{CR2}) is automatically
satisfied, so that locally $CR$-symmetric spaces are $CR$ manifolds.

\smallskip
It is proved in \cite{KZ} that a CR-symmetric space $M$ is {\it
$CR$-homogeneous}: in fact the subgroup of the automorphism
group $Aut_{CR}(M)$ generated by the symmetries
acts transitively. In particular, every $CR$-symmetric
space $M$ is a {\em real analytic} $CR$ manifold.

\bigskip
From now on we specialize to strongly pseudoconvex $CR$ manifolds
of hypersurface type and discuss $CR$-symmetric Webster metrics.
We remark that for a pseudohermitian manifold ${\cal
D}_\infty(x)=T_xM$, so that a $CR$-symmetry at a point $x$ is
characterized by the condition $(ds)_x=-Id$ on $H_xM$.

\begin{lem}
\label{lem}
Let $(M,HM,J,\eta)$ be a pseudohermitian
manifold.
Denote by  $g=g_\eta$ the Webster metric
associated to $\eta$ and by $(\varphi,\xi,\eta,g)$
the corresponding contact metric structure. Let $x\in M$
and assume that
 $\sigma_x:U\to U$ is a local $CR$-symmetry at $x$
defined on an open neighbourhood of $x$.
Then $\sigma_x$ is local
automorphism of $(\varphi,\xi,\eta,g)$.
\end{lem}

\begin{proof}
According to \cite{KZ}, Remark 3.4, we see that
the differential of $\sigma_x$ at $x$ is given by
\begin{equation}
\label{dsigma}
(d\sigma_x)_x=-Id + 2\eta_x\otimes\xi_x.
\end{equation}
Thus $(d\sigma_x)_x(\xi_x)=\xi_x$ which
implies that $(\sigma_x)_*\xi=\xi$ because $\sigma_x$
is a $CR$-isometry. It also follows that $(\sigma_x)^*\eta=\eta$
because $\eta$ is dual to $\xi$ with respect to $g$.
Since $\sigma_x$ is a $CR$ map
it follows immediately that it also preserves
the tensor field $\varphi$.
\end{proof}

\medskip
At this point we get the following characterization of
$CR$-symmetric Webster metrics:

\begin{theo}
\label{equivalent}
Let $(M,HM,J,\eta)$ be a pseudohermitian
manifold. Assume that the Webster metric $g_\eta$
is not Sasakian. The following conditions
are equivalent:

\smallskip
a) The Webster metric $g_\eta$ is locally $CR$-symmetric.

\smallskip
b) The underlying contact metric structure
satisfies the $(k,\mu)$-nullity condition.
\end{theo}

\begin{proof}
$a)\Rightarrow b)$.
It suffices to prove that the contact metric
structure is $\eta$-parallel.
Let $X,Y,Z\in{\cal D}$. We need to prove that $g((\nabla_Xh)Y,Z)=0.$
Fix a point $x\in M$ and consider a local $CR$-symmetry
$\sigma_x$ at $x$. According to the Lemma, $\sigma$ preserves
the tensor field $h$ and also its covariant derivative.
Hence at $x$ we obtain
$$g_x((\nabla_{X_x}h)Y_x,Z_x)=
g_x((\nabla_{d\sigma_x(X_x)}h)d\sigma_x(Y_x),d\sigma_x(Z_x))=
-g_x((\nabla_{X_x}h)Y_x,Z_x)$$
and the assertion follows.

\smallskip
$b)\Rightarrow a)$.
In \cite{BO1} the following tensor field $T$ is considered:
$$T_XY=(g(\varphi X,Y)+g(\varphi hX,Y))\xi-
\eta(Y)(\varphi X+\varphi hX)-\frac{\mu}{2}\eta(X)\varphi Y.$$
$T$ is a homogeneous structure for the contact metric
structure $(\varphi,\xi,\eta,g)$, i.e.
$$\bar\nabla g=\bar\nabla R=\bar\nabla T=0$$
$$\bar\nabla\xi=\bar\nabla\eta=\bar\nabla\varphi=0$$
where $\bar\nabla=\nabla-T$, $\nabla$ being the
Levi-Civita connection and $R$ its curvature tensor.
We also remark
that $\bar\nabla\bar R=\bar\nabla\bar T=0$, where $\bar R$
is the curvature tensor and $\bar T$ is the torsion of $\bar\nabla$.
Fix a point $x\in M$. From the expression of $T$
we see that $T_x$ is preserved by the linear
transformation $L:=-Id+2\eta_x\otimes\xi_x$ of $T_xM$.
Moreover, from the expression (\ref{Rkmu}) of the curvature tensor
$R$, it is straightforward to verify that $L$ also
preserves $R_x$.
This in turn implies that $L$
preserves $\bar R_x$.
Hence by a standard argument (cf. \cite{KN} p. 261)
there exists an affine transformation $\sigma:U\to U$ with respect to
$\bar\nabla$, defined on an
open neighbourhood of $x$, such that $(d\sigma)_x=L.$
From the parallelism of the structure tensors $(\varphi,\xi,\eta,g)$,
it follows that $\sigma$ is actually a $CR$-isometry, and
thus a local $CR$-symmetry at $x$.
\end{proof}

\medskip
To end this section, we shall prove that for Sasakian
manifolds, local $CR$-symmetry is actually equivalent to
a similar concept in literature, namely {\em locally $\varphi$-symmetric}
contact metric structure (cf. \cite{BLAIR+},\cite{BO1}). The latter
is defined by the requirement that the characteristic reflections,
i.e. the reflections with respect to the integral curves of $\xi$,
be local isometries. A (global) {\em Sasakian $\varphi$-symmetric space}
is a Sasakian locally $\varphi$-symmetric space whose characteristic
reflections are globally defined and $\xi$ generates
a global one-parameter group
of automorphisms of the contact structure \cite{TA}.

\begin{prop}
Let $(M,HM,J,\eta)$ be a pseudohermitian
manifold. Assume that the Webster metric $g_\eta$
is Sasakian. The following conditions
are equivalent:

\smallskip
a) $g_\eta$ is locally (globally) $CR$-symmetric.

\smallskip
b) $M$ is a locally (globally) $\varphi$-symmetric space.

\end{prop}

\begin{proof}
We treat the local statement first.

a) $\Rightarrow$ b). Since the metric $g_\eta$ is Sasakian, it is
known that the geodesic reflection $s_x$ at a point $x\in M$ with
respect to the integral curve of $\xi$ through $x$ is given, on a
normal neighbourhood $U$ of $x$ by
$$s_x=\exp_x\circ L\circ\exp^{-1}_x$$
where $L=-Id+2\eta_x\otimes\xi_x$ \cite{BV}. By (\ref{dsigma}) in
Lemma \ref{lem} the $CR$-symmetry $\sigma$ at $x$ coincides with
$s_x$ on a suitable $U'\subset U$. Hence $s_x$ is a local isometry.
This means that $M$ is locally $\varphi$-symmetric.

b) $\Rightarrow$ a) Under the assumption b), it is proved in
\cite{BV} that $s_x$ is a local automorphism of the contact metric
structure, in particular it is a $CR$-isometry and $(ds)_x=-Id$ on
$H_xM$, so that $s_x$ is a $CR$-symmetry at $x$.

Finally, as for the global statement, we remark that
if $g_\eta$ is globally $CR$-symmetric, $M$ is $CR$-homogeneous
and Riemannian homogeneous,
hence $\xi$ is complete, being a Killing field.
\end{proof}

\section{The Bochner type tensor of a $CR$-symmetric manifold}
\label{sectboch}

Let $(M,HM,J)$ be a strongly pseudoconvex $CR$ manifold
having $CR$-dimen\-sion $n\ge 2$
and let $\eta$ and $\eta'$ be two pseudohermitian
structures, with subordinate contact metric structures
$(\varphi,\xi,\eta,g)$ and
$(\varphi',\xi',\eta',g')$.
As proved in \cite{SakTak1}, these structures are
related by
\begin{equation}\label{change}
\quad\eta'=e^{2\mu}\eta ,\quad \xi'= e^{-2\mu}(\xi+Q),\quad
\varphi'=\varphi+\eta\otimes P,
\end{equation}
\[\quad g'(X,Y)=e^{2\mu}g(X,Y)\quad\forall\, X,Y\in{\cal D}
\]
where $\mu$ is a ${\cal C}^\infty$-function, $P\in{\cal D}$ is
defined by $g(P,X)=d\mu(X)$ for $X\in{\cal D}$ and $Q=JP$. In
\cite{SakTak2} the authors derive a pseudoconformal invariant on the
$CR$ manifold, that is an invariant of the change (\ref{change}),
called \emph{Bochner curvature tensor}.
A more general treatment is
given in \cite{TanB}, where almost $CR$ manifolds are allowed.
There it is proved that for $CR$ manifolds this tensor coincides with the
Chern-Moser-Tanaka invariant (\cite{CM},\cite{tanhyp}).
The definition of the Bochner curvature
tensor involves the curvature of the canonical connection
$\tilde\nabla$, as
described in the following.

As usual, the Ricci tensor field $s$ of $\tilde\nabla$ is defined by
\[s(X,Y)=\mathrm{tr}(V\rightarrow \tilde R(V,X)Y)\]
for any $X,Y\in{\frak X}(M)$. One can define another Ricci tensor
field $k$ by
\[k(X,Y)=\frac{1}{2}\,\mathrm{tr}(\varphi \tilde R(X,\varphi Y))\]
for any $X,Y\in{\frak X}(M)$. Both $s$ and $k$ are symmetric when
restricted to $\cal D$ and they satisfy
\begin{equation}\label{k=s+}
k(X,Y)=s(X,Y)+2(n-1)g(FJX,Y) \end{equation} for any $X,Y\in{\cal
D}$. We shall also denote by $\rho$ the {\em  Webster scalar curvature}
which is defined by
\[\rho=\mathrm{tr}(s).\] The expression of the Bochner curvature tensor
also involves the following tensors $l$ and $m$ defined by
\[l(X,Y)={}-\frac{1}{2(n+2)}\,k(X,Y)+\frac{1}{8(n+1)(n+2)}\,\rho\, g(X,Y)\]
\[m(X,Y)=l(JX,Y)\]
for any $X,Y\in{\cal D}$, and the tensors $L$ and $M$ such that
\[g(LX,Y)=l(X,Y),\quad g(MX,Y)=m(X,Y),\]
which satisfy $LJ=JL=M$. After this, the Bochner curvature tensor is
defined by
\[B=B_0+B_1,\]
where, for any $X,Y,Z\in{\cal D}$
\begin{eqnarray}\label{B0}
B_0(X,Y)Z\!\!&=&\!\!\tilde R(X,Y)Z-2\{m(X,Y)JZ+g(JX,Y)MZ\}\\{}
&&\!\!{}+l(Y,Z)X-l(X,Z)Y+m(Y,Z)JX-m(X,Z)JY\nonumber\\
&&\!\!{}+g(Y,Z)LX-g(X,Z)LY+g(JY,Z)MX-g(JX,Z)MY,\nonumber
\end{eqnarray}
\begin{equation}\label{B1}
B_1(X,Y)Z=\frac{1}{2}\{\tilde R(JX,JY)Z-\tilde R(X,Y)Z\}.
\end{equation}

\begin{rem}
\emph{In \cite{SakTak1} and \cite{SakTak2} the authors actually consider the
canonical connection associated to the structure
$(\bar\varphi,\bar\xi,\bar\eta,\bar g)$ such that
\[\bar\varphi=-\varphi,\quad\bar\xi=\frac{1}{2}\,\xi,\quad\bar\eta=2\eta,
\quad\bar g=4g.\]
It can be
easily seen that the connections associated to $(\bar\varphi,\bar\xi,\bar\eta,\bar g)$ and to $(\varphi,\xi,\eta,g)$
through conditions in (\ref{canonical}) coincide.
Since our
computations involve $(k,\mu)$-spaces, we prefer to express the
Bochner curvature tensor in terms of $(\varphi,\xi,\eta,g)$.}
\end{rem}

\begin{lem}
Let $(M,HM,J,\eta)$ be a pseudohermitian manifold
with associated contact metric structure
$(\varphi,\xi,\eta,g)$. Let $\nabla$ be the Levi-Civita connection
of $g$ and $\tilde\nabla$ the canonical connection. Then the
corresponding  curvature
tensors $R$ and $\tilde R$ are related by:
\begin{eqnarray}\label{Rcanonic}
\tilde R(X,Y)Z&=&R(X,Y)Z+g(FY-\varphi Y,Z)(FX-\varphi X)\\
&&{}-g(FX-\varphi X,Z)(FY-\varphi Y)-2g(\varphi X,Y)\varphi
Z\nonumber\\
&&{}+g((\tilde\nabla_XF)Y-(\tilde\nabla_YF)X,Z)\xi\nonumber
\end{eqnarray}
for any $X,Y,Z\in{\cal D}$. Denoting by $Ric$ the Ricci tensor of
$g$, then the Ricci tensor $s$ satisfies:
\begin{equation}\label{Ricci-canonic}
s(X,Y)=Ric(X,Y)-g(R(X,\xi)\xi,Y)-g(F^2X,Y)+3g(X,Y) \end{equation}
for any $X,Y\in{\cal D}$.
Finally, denoting by $\tau$ the scalar curvature of $g$,
the Webster scalar curvature is given by:
\begin{equation}
\label{rho}
\rho=\tau-2Ric(\xi,\xi)-\mathrm{tr}(F^2)+6n.
\end{equation}
\end{lem}

\begin{proof}
Straightforward computations using (\ref{H}) and the
parallelism of the structure tensors with respect to $\tilde\nabla$.
\end{proof}

\medskip
\begin{theo}\label{space form}
A pseudohermitian space form of $CR$-dimension
$n\ge 2$ is a spherical $CR$ manifold.
\end{theo}

\begin{proof}
Consider a pseudohermitian space form $(M,HM,J,\eta)$.
We need to prove that the Bochner curvature tensor vanishes.
Since the pseudoholomorphic sectional curvature is a constant $c$,
by Prop. 5.2 in \cite{Cho},
using (\ref{Rcanonic}) we obtain the following formula
for the curvature tensor of the canonical
connection $\tilde\nabla$:
\begin{align*}
\tilde R (X,Y)Z &= \frac{c}{4}\,\{g(Y,Z)X-g(X,Z)Y\\&\quad{}+g(\varphi
Y,Z)\varphi
X-g(\varphi X,Z)\varphi Y+2g(X,\varphi Y)\varphi Z\}\\
&\quad{}+g(hY,Z)X-g(hX,Z)Y+g(\varphi hY,Z)\varphi X-g(\varphi
hX,Z)\varphi Y\\
&\quad{}+g(Y,Z)hX-g(X,Z)hY+g(\varphi Y,Z)\varphi hX-g(\varphi
X,Z)\varphi hY
\end{align*}
for any $X,Y,Z\in\cal{D}$. Using (\ref{B1}), a straightforward
computation shows that
\begin{align*}
B_1(X,Y)Z &= g(hX,Z)Y-g(hY,Z)X+g(\varphi hX,Z)\varphi Y-g(\varphi hY,Z)\varphi X \\
&\quad{}+g(X,Z)hY-g(Y,Z)hX+g(\varphi X,Z)\varphi hY-g(\varphi
Y,Z)\varphi hX.
\end{align*}
Now, taking $X,Y\in\cal D$, observing that $g(\tilde R(\xi,X)Y,\xi)=0$
and $\mathrm{tr}(h)=\mathrm{tr}
(h\varphi)=0$, for the Ricci tensor field $s$ we obtain
$$s(X,Y)=\frac{c}{2}(n+1)g(X,Y)+2(n-1)g(hX,Y).$$
Applying (\ref{k=s+}), since $F=h\varphi$, for the Ricci
tensor field $k$ we obtain the following expression:
\[k(X,Y)=\frac{c}{2}(n+1)g(X,Y).\]
Computing the Webster scalar curvature, we have
\[\rho=cn(n+1).\]
With these ingredients one can compute the tensor fields $l$, $m$, $L$
and $M$ which are given by
\[l(X,Y)=-\frac{c}{8}g(X,Y),\quad m(X,Y)=-\frac{c}{8}g(\varphi
X,Y),\] \[LX=-\frac{c}{8}X,\quad MX=-\frac{c}{8}\varphi X.\]
Applying (\ref{B0}), we have
\begin{align*}
B_0(X,Y)Z&=g(hY,Z)X-g(hX,Z)Y+g(\varphi hY,Z)\varphi X-g(\varphi
hX,Z)\varphi Y\\
&\quad{}+g(Y,Z)hX-g(X,Z)hY+g(\varphi Y,Z)\varphi hX-g(\varphi
X,Z)\varphi hY\\&\quad +\frac{c}{4}\,\{g(Y,Z)X-g(X,Z)Y+g(\varphi Y,Z)\varphi
X-g(\varphi X,Z)\varphi Y\\&\quad{}+2g(X,\varphi Y)\varphi
Z\}+\frac{c}{2}\,g(\varphi X,Y)\varphi Z\\ &\quad
{}-\frac{c}{4}\,\{g(Y,Z)X-g(X,Z)Y+g(\varphi Y,Z)\varphi X-g(\varphi
X,Z)\varphi Y\}\\
&=g(hY,Z)X-g(hX,Z)Y+g(\varphi hY,Z)\varphi X-g(\varphi
hX,Z)\varphi Y\\
&\quad{}+g(Y,Z)hX-g(X,Z)hY+g(\varphi Y,Z)\varphi hX-g(\varphi
X,Z)\varphi hY.
\end{align*}
It follows that $B=B_0+B_1=0$.
\end{proof}

\begin{theo}
\label{B=0}
Let $(M,HM,J,\eta)$ be a non Sasakian locally $CR$-symmetric pseudohermitian
manifold having $CR$-dimension $n\ge 2$. Let
$(\varphi,\xi,\eta,g)$ be the underlying contact
metric structure.
Then the Bochner
curvature tensor is given by
\begin{eqnarray}\label{bochner(k,mu)}
B(X,Y)Z &=& \frac{\rho}{4n^2(n+1)}\,(g(Y,Z)X-g(X,Z)Y\\&&{}+g(\varphi
Y,Z)\varphi X-g(\varphi X,Z)\varphi Y-2g(\varphi X,Y)\varphi Z)\nonumber\\
&&{}+\frac{\rho}{2n\,\mathrm{tr}(h^2)}\,(g(hY,Z)hX-g(hX,Z)hY\nonumber\\&&{}+g(\varphi
hY,Z)\varphi hX-g(\varphi hX,Z)\varphi hY)\nonumber
\end{eqnarray}
\noindent for any $X,Y,Z\in{\cal D}$. Moreover, the following
conditions are equivalent:

\medskip
i) $B=0$.

ii) The Boeckx invariant $I=0$.

iii) The Webster scalar curvature $\rho$ vanishes.

iv) $M$ has constant pseudoholomorphic curvature.
%

\noindent
If any of the above conditions holds, then $\tilde K=0$, but
$\tilde R\not=0$.
\end{theo}

\begin{proof}
First we compute the curvature $\tilde R$ of the canonical connection
$\tilde\nabla$. Let us consider $X,Y,Z\in{\cal D}$. We remark
that, since $M$ is a $(k,\mu)$-space we have
$R(X,Y)Z\in{\cal D}$. Hence from
(\ref{Rcanonic}) we obtain
\begin{eqnarray*}
\tilde R(X,Y)Z&=&R(X,Y)Z+g(\varphi hY+\varphi Y,Z)(\varphi hX+\varphi X)\\
&&{}-g(\varphi hX+\varphi X,Z)(\varphi hY+\varphi Y)-2g(\varphi
X,Y)\varphi Z,
\end{eqnarray*}
where we applied $F=-\varphi h$. Hence, taking into account  (\ref{Rkmu})
we get
\begin{eqnarray*}
\tilde R(X,Y)Z &=&\left(1-\frac{\mu}{2}\right)(g(Y,Z)X-g(X,Z)Y\\&&{}
+g(\varphi Y,Z)\varphi X-g(\varphi X,Z)\varphi Y-2g(\varphi
X,Y)\varphi Z)\\
&&{}+\frac{1-\frac{\mu}{2}}{1-k}(g(hY,Z)hX-g(hX,Z)hY\\&&{}
+g(\varphi
hY,Z)\varphi hX-g(\varphi hX,Z)\varphi hY)\\
&&{}+g(Y,Z)hX-g(X,Z)hY+g(hY,Z)X-g(hX,Z)Y\\
&&{}+g(\varphi hY,Z)\varphi X-g(\varphi hX,Z)\varphi
Y\\&&{}+g(\varphi Y,Z)\varphi hX-g(\varphi X,Z)\varphi hY.
\end{eqnarray*}
Computing $\tilde R(\varphi X,\varphi Y)Z$, from (\ref{B1}) we get
\begin{eqnarray*}
B_1(X,Y)Z \!\!&=&\!\!\!\! {}-g(Y,Z)hX+g(X,Z)hY+g(\varphi Y,Z)h\varphi X
-g(\varphi X,Z)h\varphi Y\\
\!\!&&\!\!\!\!{}-g(hY,Z)X+g(hX,Z)Y+g(h\varphi Y,Z)\varphi
X-g(h\varphi X,Z)\varphi Y.
\end{eqnarray*}
The Riemannian Ricci tensor is given by  (cf. \cite{BKP})
\[Ric(X,Y)=(2(n-1)-n\mu)g(X,Y)+(2(n-1)+\mu)g(hX,Y)\]
for any $X,Y\in{\cal D}$. Using (\ref{Ricci-canonic}), since
$F^2X=h^2X=(1-k)X$, we get
\[s(X,Y)=n(2-\mu)g(X,Y)+2(n-1)g(hX,Y),\]
and the Webster scalar curvature is
\begin{equation}
\label{rhokmu}
\rho=2n^2(2-\mu),
\end{equation}
which proves that $I$ vanishes if and only if $\rho=0$.
Applying (\ref{k=s+}),
\[k(X,Y)=n(2-\mu)g(X,Y).\]
With these elements, the tensors $l,m,L,M$ are given by
\[l(X,Y)=\frac{n(\mu-2)}{4(n+1)}g(X,Y),
\quad m(X,Y)=\frac{n(\mu-2)}{4(n+1)}g(\varphi X,Y)\]
\[LX=\frac{n(\mu-2)}{4(n+1)}X,
\quad MX=\frac{n(\mu-2)}{4(n+1)}\varphi X.\] Using (\ref{B0}) and the
expression for $\tilde R$, a straightforward computation shows that
\begin{eqnarray*}
B_0(X,Y)Z\!\!&=&\!\!\!\!
\frac{2-\mu}{2(n+1)}(g(Y,Z)X-g(X,Z)Y\\\!\!&&\!\!\!\!{} +g(\varphi
Y,Z)\varphi X-g(\varphi X,Z)\varphi Y-2g(\varphi
X,Y)\varphi Z)\\
\!\!&&\!\!\!\!{}+\frac{2-\mu}{2(1-k)}(g(hY,Z)hX-g(hX,Z)hY\\\!\!&&\!\!\!\!{}
+g(\varphi
hY,Z)\varphi hX-g(\varphi hX,Z)\varphi hY)\\
\!\!&&\!\!\!\!{}+g(Y,Z)hX-g(X,Z)hY-g(\varphi Y,Z)h\varphi X+g(\varphi X,Z)h\varphi Y\\
\!\!&&\!\!\!\!{}+g(hY,Z)X-g(hX,Z)Y-g(h\varphi Y,Z)\varphi
X+g(h\varphi X,Z)\varphi Y,
\end{eqnarray*}
and we get the expression (\ref{bochner(k,mu)})
for $B=B_0+B_1$, since
$\mathrm{tr}(h^2)=2n(1-k).$
The equivalence of i) and iii) is an immediate consequence.
As for the pseudoholomorphic sectional curvature, consider a holomorphic
2-plane $\sigma=<X,JX>$ where $X$ is a unit holomorphic tangent vector at some
point $x\in M$.
Using the
expression of $\tilde R$, we get
\begin{eqnarray*}
\tilde K(\sigma) &=& g(\tilde R(X,\varphi X)\varphi X,X)\\
&=& 2(2-\mu)-\frac{2-\mu}{1-k}(g(hX,X)^2+g(h\varphi X,X)^2)
\end{eqnarray*}
showing that $\tilde K$ vanishes for
$\mu=2$. Conversely, if $\tilde K$ is constant, Theorem
\ref{space form} guarantees that $B=0$.
\end{proof}

\section{The classification}
In this section we prove Theorem \ref{main} and its Corollaries.

\smallskip
\begin{proof}{\sc of Theorem \ref{main}}\
Let  $(M,HM,J,\eta)$ be a simply connected pseudohermitian
manifold which is a spherical $CR$ manifold and such that $g_\eta$
is $CR$-symmetric. If $g_\eta$ is not Sasakian, according to
Theorem \ref{B=0} $M$ is a $(k,\mu)$-space with vanishing Boeckx
invariant. It follows that $M$ is homothetic to $T_1{\mathbb H}^{n+1}$
endowed with its standard $CR$ structure, since it is known that
the Boeckx invariant of $T_1{\mathbb H}^{n+1}$ vanishes \cite{BKP}.
Next we consider the case where $g_\eta$ is
Sasakian. Then $M$ is a simply connected Sasakian
$\varphi$-symmetric space and according to \cite{JK} it is a
principal fiber bundle $\pi:M\to N$ over a simply connected
Hermitian symmetric space $(N,g_o)$ and $\pi:(M,g_\eta)\to (N,g_o)$ is
a Riemannian submersion with fibers tangent to $\xi$, which is
also a $CR$ map. Moreover, since $M$ is spherical, $N$ is
Bochner-flat. Indeed, denoting by $B_N$ the Bochner tensor of $N$,
by a result of S. Webster \cite{Web} already quoted in the
Introduction, we have
$$\pi_*(B(X,Y)Z)=B_N(\pi_*X,\pi_*Y)\pi_*Z$$
for any $x\in M$ and $X,Y,Z\in H_xM$.
Now, according to
a result of M. Matsumoto and S. Tanno \cite{MT}, $N$ is
either a simply connected K\"ahler space form or is
isometric to
a product $N^k(c)\times N^{n-k}(-c)$, $c>0$, of two
simply connected K\"ahler space forms with holomorphic curvatures
respectively $c$ and $-c$.
In the first case, $M$ is a Sasakian space form
and hence, as a pseudohermitian manifold, it is
homothetic to $\mathbb S^{2n+1}$, $H^{2n+1}$,
or $B^n\times\mathbb R$. In the last case,
up to a homothetic change of the metric $g_o$, we can assume $c=1$,
yielding a homothetic change of the pseudohermitian
structure $\eta$ of $M$ which turns $M$ into a Sasakian
manifold equivalent to
the Sasakian $\varphi$-symmetric space $P^n_k$ over
${\mathbb C\mathbb P}^k\times{\mathbb C\mathbb H}^{n-k}$.
\end{proof}

\medskip

The proof of Corollary \ref{cor1} is immediate.

\medskip
\begin{proof}{\sc of Corollary \ref{cor2}}
According to Theorem \ref{B=0},
the assumption on $\tilde K$ forces
$g_\eta$ to be Sasakian, hence $M$ is a Sasakian
$\varphi$-simmetric space. Thus the simply connected
covering $\tilde M$ is also a $\varphi$-symmetric
space which is locally equivalent to $M$ as a pseudohermitian
manifold. In particular, $\tilde M$ is a spherical
$CR$-symmetric space with positive pseudoholomorphic curvature at some
point, which is a principal fiber bundle $\pi:\tilde M\to N$
over a Hermitian symmetric space $N$.
Now, comparing with the classification in Theorem \ref{main},
$\tilde M$ cannot be homothetic to any of
the spaces $P^n_k$. Indeed,
observe that at each point $x$ of $P^n_k$ we have
$\tilde K(\sigma)=0$ for some
holomorphic 2-plane $\sigma$. Indeed,
choose a holomorphic
2-plane $\sigma'$ of $\mathbb C\mathbb P^k\times\mathbb C\mathbb H^{n-k}$ at $\pi(x)$
with vanishing holomorphic curvature. Such a $\sigma'$ exists
since $\mathbb C\mathbb P^k$ and $\mathbb C\mathbb H^{n-k}$
have opposite holomorphic curvatures.
Now take $\sigma$ such
that $\pi_*(\sigma)=\sigma'$.
For the other models in the classification except for
$\mathbb S^{2n+1}$, it is known that at each point $\tilde K\le 0$.
Hence $\tilde M$ is compact and the assertion follows.
\end{proof}

\section{Examples: $CR$ geometry of tangent sphere bundles}

Let $(M,g)$ be a Riemannian manifold of
dimension $n\ge 3$. For each $r>0$ we shall
denote by $T_rM$ the tangent sphere bundle of radius $r$,
which is the hypersurface of the tangent bundle $TM$
defined by
$$T_rM=\{(x,u)\in TM|g_x(u,u)=r^2\}.$$
Here and in the following we consider a point of the
tangent bundle $TM$ as a pair $(x,u)$ with $x\in M$
and $u\in T_xM$.
Let $\pi:TM\to M$ be the canonical projection such that
$\pi(x,u)=x$.
For each smooth vector field $X\in\frak X(M)$ we shall
denote by $X^V$ its {\em vertical lift}
to $TM$ and by $X^H$ its {\em horizontal lift} with respect to
the Levi-Civita connection $D$ of $(M,g)$. For details,
our standard reference is \cite{BLAIR+}, Ch. 9.
If $t=(x,u)$ is a fixed point of $TM$ and $X\in T_xM$, we shall
also denote by $X_t^H\in T_tTM$ its horizontal lift and
by $X_t^V\in  T_tTM$ its vertical lift.
Then at each point $t=(x,u)$ of $T_rM$ the tangent
space to $T_rM$ at $t$ is given by
$$T_t(T_rM)=\{X_t^H+Y_t^V|X,Y\in T_xM,\,\,g_x(Y,u)=0\}.$$

\medskip
Let $\lambda\not=0$ be a fixed real number.
One can define an almost complex structure $J_\lambda:TTM\to TTM$
by
\begin{equation}
\label{Jlambda}
J_\lambda(X^H)=\lambda X^V,\,J_\lambda(X^V)=-\frac{1}{\lambda}X^H.
\end{equation}
Since $T_rM$ is a real hypersurface of $TM$, it inherits
canonically
a partial complex structure $(H(T_rM),J_\lambda)$ from $J_\lambda$.
The  holomorphic tangent bundle $H(T_rM)$ can be described as
follows.
At a fixed point $t=(x,u)\in T_rM$ we have:
\begin{equation}
\label{HTrM}
H_t(T_rM)=\{X_t^H+Y_t^V|X,Y\in T_xM,\,\,g_x(X,u)=g_x(Y,u)=0\}.
\end{equation}
We define a global horizontal vector field $\xi\in\frak X(T_rM)$ by
\begin{equation}
\label{xi}
\xi_t=\frac{2}{\lambda}u_t^H,\quad t=(x,u).
\end{equation}
Denoting by $\frak U$ the canonical vertical vector
field of $TM$ (cf. \cite{BLAIR+}, p. 142 or \cite{KOWSEK} p. 210)
we have that
$$\xi=-2J_\lambda\frak U.$$
We recall that the local expression of $\frak U$ in
a coordinate system $(x^i,v^i)$ of $TM$ induced
by a local chart $(U,x^1,\dots, x^n)$ of $M$ is
$$\frak U=v^i\frac{\partial}{\partial v^i}.$$
Then $\xi$ is everywhere transverse to the holomorphic
tangent bundle $H(T_rM)$.

\begin{theo}
For each $r>0$ and $\lambda>0$, $(T_rM, H(T_rM),J_\lambda)$ is
a strictly pseudoconvex almost $CR$ manifold.
Moreover, if $(M,g)$ has constant curvature then $T_rM$
is locally $CR$-symmetric with respect to the Webster metric
$g_{\eta_\lambda}$ where $\eta_\lambda$ is the pseudohermitian structure
such that $\eta_\lambda(\xi)=1.$ If, in addition, $M$ is simply connected
and complete, $g_{\eta_\lambda}$ is globally $CR$-symmetric.
\end{theo}

\begin{proof}
Define the 1-form $\eta_\lambda$ on $T_rM$ by
$$\eta_\lambda(H(T_rM))=0,\,\eta_\lambda(\xi)=1.$$
First we shall prove that the partial complex structure
$(H(T_rM),J_\lambda)$ satisfies (\ref{CR1}) and that
the Levi form associated to $\eta_\lambda$
is positive definite at each point $t=(x,u)\in T_rM$.
In order to simplify the notation, in the
following we shall denote simply by $J$ both
the almost complex structure $J_\lambda$ on $TM$ and
the partial complex structure induced on $T_rM$.
We shall also denote $\eta_\lambda$ by $\eta$.

Hence we shall verify that, for each $t=(x,u)\in T_rM$
\begin{equation}
\label{claim1}
d\eta(Z,W)=d\eta(JZ,JW),\quad \,\,d\eta(JZ,Z)>0
\end{equation}
where $Z,W\in H_t(T_rM).$
To this aim, we shall use the fact that, according to (\ref{HTrM}),
$H_t(T_rM)$ is spanned by vectors of the form $X_t^H$
and $X^V_t$ where $X\in\frak X(M)$ is such
that $X_x$ is orthogonal to $u$
with respect to $g$.
We remark that $X_t^V$ can be extended to a global
section $X^t$ of $H(T_rM)$ defined as follows.
Let $g^S$ be the Sasaki metric on $TM$ naturally
constructed from $g$ (cf. e.g. \cite{BLAIR+} or \cite{KOWSEK}).
Observe that $\xi$ is orthogonal to $H(T_rM)$
with respect to the Riemannian metric induced by $g^S$
on $T_rM$, which will be denoted by the same symbol.
Then we set
$$X^t:=X^V-\frac{1}{r^2}g^S(X^V,\frak U)\frak U.$$
The vector field $X^t$ is the {\em tangential lift}
of $X$ as defined e.g. in \cite{KOWSEK}, p. 211.

We also remark that $X_t^H$
can be extended to a global
section $X^0$ of $H(T_rM)$ defined as follows:
\begin{equation}
\label{hollift}
X^0:=X^H-\frac{\lambda}{2r^2}g^S(X^V,\frak U)\xi.
\end{equation}
Actually according to (\ref{Jlambda}) we have
\begin{equation}
\label{JonTrM}
JX^0=\lambda X^t,\,\,JX^t=-\frac{1}{\lambda}X^0.
\end{equation}

After this, we first compute $d\eta(JX_t^H,Y_t^H)$ where
$X,Y\in\frak X(M)$ with $g_x(X_x,u)=g_x(Y_x,u)=0$.
According to (\ref{JonTrM}), we have
\begin{eqnarray*}
2d\eta(JX^0,Y^0)=-\lambda\eta[X^t,Y^0].
\end{eqnarray*}
Moreover, taking into account the formula (cf. \cite{BLAIR+}, p. 138)
\begin{equation}
\label{HV}
[X^H,Y^V]=(D_XY)^V
\end{equation}
we obtain, evaluating at the point $t$:
$$\eta[X^t,Y^0](t)=\eta\left(-\frac{\lambda}{2r^2}X^Vg^S(Y^V,\frak U)\xi
\right)(t).$$
On the other hand, it is readily verified that at $t$
the function $X^Vg^S(Y^V,\frak U)$ takes the value $g_x(X,Y).$
Hence
\begin{equation}
\label{step}
\eta[X^t,Y^0](t)=-\frac{\lambda}{2r^2}g_x(X,Y)
\end{equation}
and we conclude that
\begin{equation}
\label{Leviform1}
d\eta(JX_t^H,Y_t^H)=\frac{\lambda^2}{4r^2}g_x(X,Y).
\end{equation}
Next we compute $d\eta(JX^V_t,Y_t^V)$. Using (\ref{JonTrM}) again we have:
$$
2d\eta(JX^t,Y^t)=\frac{1}{\lambda}\eta[X^0,Y^t].$$
Hence, evaluating at $t$ and taking into account (\ref{step}) we get
\begin{equation}
\label{Leviform2}
d\eta(JX^V_t,Y_t^V)=\frac{1}{4r^2}g_x(X,Y).
\end{equation}
Next observe that from
$$2d\eta(JX^0,Y^t)=-\lambda\eta[X^t,Y^t]$$
using $[X^V,Y^V]=0$, we obtain
\begin{equation}
\label{Leviform3}
d\eta(JX_t^H,Y^V_t)=0.
\end{equation}
Thus, taking into account (\ref{HTrM}), equations (\ref{Leviform1}),
(\ref{Leviform2}) and (\ref{Leviform3}) yield (\ref{claim1}).
The first assertion is proved.

Now, suppose $(M,g)$ has
constant curvature $K$. We shall compute first the expression
of the Webster metric $g_\eta$.
First of all we claim that $\xi$ is actually the Reeb
vector field of $\eta$ i.e.
$d\eta(Z,\xi)=0$ for any $Z\in\frak X(T_rM).$
To justify this, it suffices to verify that
$$d\eta(X^0,\xi)=d\eta(X^t,\xi)=0$$
for any $X\in\frak X(M)$, or equivalently
$\eta[X^0,\xi]=\eta[X^t,\xi]=0,$
which in turn is equivalent to
\begin{equation}
\label{claim2}
g^S([X^0,\xi],\xi)=g^S([X^t,\xi],\xi)=0.
\end{equation}
Indeed, we have
$$[X^0,\xi]=[X^H,\xi]+\frac{\lambda}{2r^2}\xi g^S(X^V,\frak U)\xi.$$
Now fix a point $t=(x,u)$ and consider a coordinate neighbourood $(x^i,v^i)$
around $t$;
then we have the local expression $\xi=\frac{2}{\lambda}v^k(\partial_k)^H$.
Assuming $X=X^i\partial_i$, $u=u^k({\partial_k})_x$,  we
compute

\begin{eqnarray}
&&[X^H,\xi]_t=
\frac{2}{\lambda}\{X_t^H(v^k)(\partial_k)_t^H+
u^k[X,\partial_k]^H_t-(R_x(X,u)u)^V_t\}\nonumber\\
&&=\frac{2}{\lambda}\{X_t^H(v^k)(\partial_k)_t^H+u^k(D_X\partial_k)^H_t
-u^k(D_{\partial_k}X)^H_t-(R_x(X,u)u)^V_t\}\nonumber\\
&&=-\frac{2}{\lambda}\{u^k(D_{\partial_k}X)^H_t+(R_x(X,u)u)^V_t\}
\label{X^Hxi}
\end{eqnarray}
where we have used the formula for the Lie brackets of type $[X^H,Y^H]$
in \cite{BLAIR+}, p. 138. Here $R$ denotes
the curvature tensor field of $(M,g)$.
It follows that:
\begin{equation}
g^S([X^H,\xi],\xi)(t)=-\frac{2}{\lambda}u^k g^S((D_{\partial_k}X)^H,\xi)(t).
\label{*}
\end{equation}
On the other hand,
$$g^S\left(\frac{\lambda}{2r^2}\xi g^S(X^V,\frak U)\xi,\xi\right)=
\frac{2}{\lambda}\xi g^S(X^V,\frak U).$$
Now, taking into account that with respect to
the Levi Civita connection $\nabla$ of $(TM,g^S)$ it holds
$\nabla_{\partial_k^H}\frak U=0$ and that the vertical
component of $\nabla_{\partial_k^H}X^V$ is $(D_{\partial_k}X)^V$
(cf. \cite{KOWSEK}, p. 210), we obtain
$$\xi g^S(X^V,\frak U)(t)=
\frac{2}{\lambda}u^kg^S((D_{\partial_k}X)^V,\frak U)(t)=
u^kg^S((D_{\partial_k}X)^H,\xi)(t).$$
Thus compairing with (\ref{*}) we can conclude
that $g^S([X^0,\xi],\xi)=0$. The proof of $g^S([X^t,\xi],\xi)=0$ is
similar and hence omitted for brevity.

Now we see that the Webster metric $g_\eta$ is the
restriction to $T_rM$ of the $g$-natural metric on $TM$:

\begin{equation}
\label{G}
G=\frac{1}{4r^2}g^S+\frac{\lambda^2-1}{4r^2}g^v
\end{equation}
(cf. \cite{Abb}, \cite{AbbSar} for
the general theory of $g$-natural metrics on tangent
bundles). Here $g^v$ denotes the
{\em vertical lift} of $g$ determined by
$$g^v(X^H,Y^H)=g(X,Y),\,g^v(X^H,Y^V)=g^v(X^V,Y^V)=0,\quad
X,Y\in\frak X(M).$$
This follows from the formulas (\ref{Leviform1}), (\ref{Leviform2}),
and (\ref{Leviform3}) for the Levi form at a generic
point $t=(x,u)\in T_rM$, from the fact that $G(X^H_t,\xi_t)=G(X^V_t,\xi_t)=0$
for every $X\in T_xM$ with $g_x(X,u)=0$,
and finally observing that $G(\xi,\xi)=1$.

After this, we show that at each point $t=(x,u)$ there
exists a local $CR$-symmetry of $T_rM$. Since $(M,g)$
has constant curvature, there exists a local
isometry $f:U\to U$ defined on an open neighbourood of $x$
in $M$ such that
$$f(x)=x,\,df_x(u)=u,\,df_x(X)=-X \textrm{ if }g_x(X,u)=0.$$
Indeed, the linear mapping $-Id+\frac{2}{r^2}u^\flat\otimes u$
of $(T_xM,g_x)$ preserves the curvature tensor $R_x$.
Now consider the induced mapping $F=df:TU\to TU$.
We remark that, since
$f$ is an isometry, $dF$ commutes with horizontal
and vertical lifts, i.e.
\begin{equation}
\label{claim3}
dF_s(X^H_s)=(df_y(X))^H_{F(s)},\,
dF_s(X^V_s)=(df_y(X))^V_{F(s)}
\end{equation}
at each point $s=(y,v)$ of $TU$, for every $X\in T_yM$.
This implies that $F$ is both a local isometry of $(TM,G)$
and a holomorphic mapping with respect to $J$.
In particular, $F$ restricts to a local isometry
of $T_rM$ which is also a $CR$ map. Finally, $F$
is a local $CR$-symmetry at $t$, since $F(t)=(f(x),df_x(u))=t$
and using (\ref{claim3})
again, for every $Z=X^H_t+Y^V_t\in H_t(T_rM)$ we have
$$dF_t(Z)=(df_x(X))_t^H+(df_x(Y))_t^V=-Z.$$
Finally notice that $f$ can be globally defined when
$M$ is complete and simply connected, thus $F$ is also
globally defined on $T_r(M)$.
\end{proof}

\begin{theo}
\label{lambdar}
Let $(M,g)$ be a Riemannian manifold with constant
curvature $K$ and dimension $n\ge 3$.
Fix $r>0$, $\lambda>0$ and consider
the $CR$ manifold $(T_rM,H(T_rM),J_\lambda)$ as above.
Then

\medskip
a) $T_rM$ is spherical if and only if $\lambda^2+Kr^2=0$.

\medskip

b) The metric $g_{\eta_\lambda}$ is Sasakian if and only if $\lambda^2-Kr^2=0$.

\medskip
c) When $(T_rM,H(T_rM),J_\lambda,g_{\eta_\lambda})$ is not Sasakian, its Boeckx invariant is
$$I=\frac{\lambda^2+Kr^2}{|\lambda^2-Kr^2|}.$$
Hence, when $K\not=0$, each $T_rM$ admits a one-parameter family
$(H(T_rM),J_\lambda,\eta_\lambda)$ of locally (globally for a complete, simply
connected $M$) $CR$-symmetric non homothetic
pseudohermitian structures. If $K<0$, exactly one of the underlying
$CR$ structures is spherical.
\end{theo}

\begin{proof}
We begin by computing the spectrum of the operator
$h=\frac{1}{2}{\cal L}_\xi\varphi$
where $\varphi$ is the $f$-structure extending $J:=J_\lambda$ on
$T(T_rM)$ by
$\varphi(\xi)=0$. Fix a point $t=(x,u)$ and consider a holomorphic
vector of the form $Z=X^V_t$ with $g_x(X,u)=0$. We shall verify
that $X^V_t$ is an eigenvector of $h_t$.
Indeed we shall compute $2h(X^t)=[\xi,JX^t]-J[\xi,X^t]$ and then
evaluate at $t$. Here $X^t$ is a tangential lift extending $Z$
as in the preceding proof. Now observe that
\begin{equation}
\label{hX^V}
[\xi,JX^t]_t-J[\xi,X^t]_t=-\frac{1}{\lambda}[\xi,X^0]_t-J[\xi,X^t]_t=
-\frac{1}{\lambda}[\xi,X^H]_t-J[\xi,X^V]_t.
\end{equation}
In a coordinate neighbourhood $(x^i,v^i)$ around $t$
we have, using (\ref{HV}):
$$[\xi,X^V]_t=\frac{2}{\lambda}\{(D_uX)^V_t -X^H_t\}$$
whence
$$J[\xi,X^V]_t=-\frac{2}{\lambda}\left\{\frac{1}{\lambda}
(D_{u}X)^H_t +\lambda X^V_t\right\}.$$
On the other hand, since $(M,g)$ has constant curvature $K$, (\ref{X^Hxi}) yields
$$\frac{1}{\lambda}[\xi,X^H]_t=\frac{2}{\lambda^2}\{
(D_{u}X)^H_t +Kr^2X^V_t\}.$$

Thus coming back to (\ref{hX^V}), we get
\begin{equation}
\label{heigenv1}
h(X^V_t)=\frac{\lambda^2-Kr^2}{\lambda^2}X^V_t.
\end{equation}
Since $h$ anticommutes with $J$, from this it also follows that
\begin{equation}
\label{heigenv2}
h(X^H_t)=-\frac{\lambda^2-Kr^2}{\lambda^2}X^H_t.
\end{equation}
for any $X\in T_xM$. Thus according to (\ref{HTrM}) we can
conclude that the spectrum of
$h$ is $\{0,\pm \frac{\lambda^2-Kr^2}{\lambda^2}\}$. Since
$T_rM$ is a $CR$ manifold, the assertion b) follows
directly.
To prove a), we need to compute the Webster
scalar curvature of $T_rM$.
To this aim, we shall compute the scalar
curvature of $g_\eta$, where $\eta:=\eta_\lambda$.
We shall denote
by $\nabla'$ the Levi Civita connection and by
$R'$ the curvature tensor of $g_\eta$.
Recall that $g_\eta$ is the restriction of
the $g$-natural metric $G$ in (\ref{G}); in particular,
we remark that $\pi:(T_rM,g_\eta)\to (M,\frac{\lambda^2}{4r^2}g)$ is
a Riemannian submersion. By standard arguments, we see that the
fibers of $\pi$ are totally geodesic and of constant curvature 4.
Indeed, we have the formula
\begin{eqnarray*}
R'(X^V_t,Y^V_t)Z^V_t=\frac{1}{r^2}\{g_x(Y,Z)X^V_t-
g_x(X,Z)Y^V_t\}
\\
=4\{g_\eta(Y^V_t,Z^V_t)X^V_t-
g_\eta(X^V_t,Z^V_t)Y^V_t\}.
\end{eqnarray*}
Here $t=(x,u)$ and $X,Y,Z\in T_xM$ are orthogonal to $u$.
Using the Gauss equation, this formula can be
derived from the fact that
$\bar R(X^V,Y^V)Z^V=0$, which holds for the curvature
of any $g$-natural metric on $TM$ of type $G=ag^S+bg^h+cg^v$
(see \cite{AbbSar}). Now, the scalar cuvature $\tau$ of
$(T_rM,g_\eta)$ is related to the scalare
curvatures $\tau_M$, $\hat\tau$
of $(M,\frac{\lambda^2}{4r^2}g)$ and
of the fibers of $\pi$ by
\begin{equation}
\label{scal}
\tau=\tau_M+\hat\tau-||A||^2
\end{equation}
where $A$ is the O'Neill fundamental horizontal tensor field of
the submersion $\pi$
(cf. e.g. \cite{IanPasFalc}). To compute $||A||^2$,
we fix a point $t=(x,u)$ and an orthonormal
frame $\{\frac{1}{rp}u,X_1,\dots, X_n\}$ of $(T_xM,p^2g)$ where
$p:=\frac{\lambda}{2r}$. Then we can consider the orthonormal basis
$\{\frac{1}{rp}u^H,(X_1)^H_t,\dots,(X_n)^H_t,\lambda (X_1)^V_t,\dots,\lambda
(X_n)^V_t\}$
of $(T_t(T_rM),g_\eta)$.
Now we take into account the formula
$$\nabla'_{X^H_t}Y^t=\frac{1}{2\lambda^2}(R_x(u,Y)X)^H_t+
(D_XY)^t_t$$
which can be derived from the formula of Gauss
for $T_rM$ and
the expression of the Levi-Civita connection
of $(TM,G)$  (cf. \cite{AbbSar}).
This formula yields:
$$A_{X_i^H}X_j^V=\frac{K}{2\lambda^2p^2}\delta^i_j\,u^H,\,
A_{u^H}X^V_i=-\frac{Kr^2}{2\lambda^2}X_i^H.
$$
Thus
$$||A||^2=\lambda^2\sum_{ij}g_\eta(A_{X_i^H}X_j^V,A_{X_i^H}X_j^V)+
\frac{\lambda^2}{r^2p^2}
\sum_{i}g_\eta(A_{u^H}X_i^V,A_{u^H}X_i^V)=2n\frac{K^2r^4}{\lambda^4}.$$
Hence using (\ref{scal}) we obtain
$$\tau=4n(n+1)\frac{Kr^2}{\lambda^2}+4n(n-1)-2n\frac{K^2r^4}{\lambda^4}$$
which in turn yields the following formula for the Webster scalar
curvature
$$\rho=4n^2\left(1+\frac{Kr^2}{\lambda^2}\right).$$
Hence, since the metric $g_\eta$ is locally $CR$-symmetric, assertion a)
follows from Theorem \ref{B=0}. The determination of the
Boeckx invariant is a immediate consequence of (\ref{heigenv1})
and (\ref{rhokmu}).
\end{proof}

\bigskip
\begin{proof}{\sc of
Corollary \ref{K><0}} Let $(M,g)$ be
a Riemannian manifold with constant curvature $K$ and dimension
$n\ge 3$.
The assertions a) and b) and the uniqueness assertions
are direct consequences of a) and b)
of Theorem \ref{lambdar}, setting $\lambda=1$.
To prove c), first we remark that when $K=0$  the $(T_rM,H(T_rM),J,g_\eta)$
are all non Sasakian and that the Boeckx invariant
actually does not depend on $r$, namely $I=1$.
Vice versa, assuming that the pseudohermitian manifolds
$T_rM$ are all locally homotethic, we see from a) and b) and the uniqueness
assertions that both $K>0$ and $K<0$ must be excluded.
\end{proof}

\bigskip
\bigskip
{\small
\begin{flushleft}
\textsc{Dipartimento di Matematica-Universit\`a di Bari\\
Via E. Orabona 4\\70125 BARI ITALY}\\ \textit{E-mail addresses}:
dileo@dm.uniba.it,\,lotta@dm.uniba.it
\end{flushleft}
}
\end{document}